\documentclass[12pt,reqno]{amsart}

\usepackage{setspace}

\author{Mikhail G. Katz}\address{M. Katz, Department of Mathematics,
Bar Ilan University, Ramat Gan 52900
Israel}\email{katzmik@macs.biu.ac.il}

\author{Semen S. Kutateladze}\address{S. Kutateladze, Sobolev
Institute of Mathematics, Novosibirsk State University,
Russia}\email{sskut@math.nsc.ru}

\begin{document}

%\doublespacing

\thispagestyle{empty}

\title{Edward Nelson (1932-2014)}

\maketitle

On September 10, 2014 we lost Edward Nelson, one of the most original
and controversial mathematical thinkers of our time.

Nelson was born on May 4, 1932 in Decatur, Georgia.  His father spoke
fluent Russian, having spent some time in St.~Petersburg in connection
with issues related to prisoners of war (Nelson 1995 \cite{Ne95}).
Nelson was educated at the University of Chicago.  In 1953 he was
awarded a master's degree, and in 1955 he earned a doctorate under the
supervision of Irving Segal, an expert in mathematical methods of
quantum mechanics and functional analysis and one of the founders of
the theory of C*-algebras.  Not surprisingly, Nelson's early papers
were devoted to quantum mechanics.  In 1995, Nelson received the
Steele prize awarded by the American Mathematical Society, for his
1963 and 1977 papers on applications of probability theory to quantum
fluctuations.

After graduating, Nelson spent three years at Princeton's Institute
for Advanced Study, and from 1959 until his retirement in 2013, he
worked at Princeton University where he became a Professor in 1964.
In his memory the university flag was raised at half-mast for three
days.  In 1975 Nelson was elected to the American Academy of Arts and
Sciences, and in 1997, to the National Academy of Sciences of the USA.
Princeton University Press published the following six monographs by
Nelson:

\begin{enumerate}
\item
Dynamical theories of Brownian motion. Princeton University Press,
Princeton, N.J. 1967.
\item
Tensor analysis.  Princeton University Press, Princeton, N.J. 1967.
\item
Topics in dynamics. I: Flows. Mathematical Notes. Princeton University
Press, Princeton, N.J.; University of Tokyo Press, Tokyo 1969.
\item
Quantum fluctuations.  Princeton Series in Physics.  Princeton
University Press, Princeton, N.J., 1985.
\item
Predicative arithmetic.  Mathematical Notes, 32.  Princeton University
Press, Princeton, N.J. 1986.
\item
Radically elementary probability theory.  Annals of Mathematics
Studies, 117.  Princeton University Press, Princeton, N.J. 1987.
\end{enumerate}
All of these monographs are freely available online at

https://web.math.princeton.edu/$\sim$nelson/books.html

Already the short list above indicates the wide range of his interests
and results related to probability theory, quantum mechanics, the
theory of dynamical systems, and mathematical logic.  The volume
\cite{Fa06} edited by Faris contains helpful surveys of Nelson's
mathematical work in varied areas as well as a full bibliography.

Exceptional fame came to Nelson in the wake of his developing a
fundamentally new, syntactic approach to Abraham Robinson's
nonstandard analysis \cite{Ro66}, based on the internal set theory
(IST) Nelson proposed in \cite{Ne77}.  A related framework was
developed independently by K.~Hrbacek \cite{Hr78}.  Detailed
presentations of IST may be found in \cite{GKK} and~\cite{KR04}.

Nelson rejected the prevailing views to the effect that nonstandard
analysis operates with some fictional elements that extend the
standard world of mathematical entities.  In his approach, nonstandard
objects inhabit the realm of the most ordinary mathematical objects.
Nelson emphasized the creative syntactic contribution of the new
approach in the following terms: ``Really new in nonstandard analysis
are not theorems but the notions, i.e., external predicates.''
(Nelson 1988 \cite[p.~134]{Ne88})

A major event of the mathematical scene was Nelson's book
\emph{Radically elementary probability theory} (Nelson 1987
\cite{Ne87}) which contains a modern exposition of von Mises frequency
approach based on internal set theory IST.  In the preface to this
book, which was translated into French and Russian, Nelson writes:
\begin{quote}

Let me take this occasion to express a strong hope. This book is a
bare beginning, and it will not succeed unless others carry on the
job.  There is a tendency in each profession to veil the mysteries in
a language inaccessible to outsiders. Let us oppose this. Probability
theory is used by many who are not mathematicians, and it is essential
to develop it in the plainest and most understandable way possible.
My hope is that you, readers of this book, will build a pavilion of
probability open to all who search for understanding. (From the
preface to the Russian edition.)
\end{quote}

To elaborate on Nelson's IST approach to infinitesimals in a
non-technical way, note that the general mathematical public often
takes the Zermelo-Fraenkel theory with the Axiom of Choice (ZFC) to be
\emph{the} foundation of mathematics.  How much ontological reality
one assigns to this varies from person to person.  Some mathematicians
distance themselves from any kind of ontological endorsement, which is
a formalist position in line with Robinson's.  On the other hand, many
do assume that the ultimate test of truth and/or verifiability is in
the context of ZFC, so that in this sense ZFC still is \emph{the
foundation}, though it is unclear whether Robinson himself would have
subscribed to such a view.  

Much of the mathematical public attaches more significance to ZFC than
merely a formalist acceptance of it as the ultimate test of
provability, and tends to accept the existence of sets (wherever they
are to be found exactly but that is a separate question), including
infinite sets, and particularly the existence of the set of
\emph{real} numbers which are often assumed to be built into the
nature of \emph{reality} itself, as it were.  When such views are
challenged what one often hears in response is an indignant monologue
concerning the categoricity of the real numbers, etc.  As G.~F.~Lawler
put it in his ``Comments on Edward Nelson's `Internal Set Theory: a
new approach to nonstandard analysis'{}'':
\begin{quote}
Clearly, the real numbers exist and have these properties. Indeed,
many courses in elementary analysis choose not to construct the reals
but rather to take the existence of an ordered field as given. This is
reasonable: we are implicitly assuming such an object exists,
otherwise, why are we studying it?  (Lawler 2011 \cite{La11})

\end{quote}
In such a situation infinitesimals may not thrive: the real numbers
are thought truly to \emph{exist}, infinitesimals \emph{not}.  This is
where Nelson comes in with his syntactic novelty and declares: Guess
what, we can find infinitesimals \emph{within the real numbers
themselves} if we are only willing to enrich the language we speak!
The perceived ontological differences between the real numbers and
infinitesimals are therefore seen to be merely a function of the
technical choices, including syntactic limitations, imposed by Cantor,
Dedekind, Weierstrass, and their followers in establishing the
foundations of analysis purged of infinitesimals.  This and related
issues are explored in \cite{Ba13} and \cite{Ba15}.

Nelson was a unique mathematician, whose creative interests shifted
from the classical areas of mathematical physics and functional
analysis toward mathematical logic and computer science under the
influence of technological progress of modern computation.  Nelson was
the first webmaster of the Mathematics Department at Princeton
University.  Free in his \emph{Ars Inveniendi}, he advocated freedom
of knowledge and information.  His motto was the paradoxical thesis
``Intellectual property is an oxymoron."

His memory will endure in our hearts.

\end{document}